\documentclass[12pt]{amsart}
\usepackage{amsmath,amscd,latexsym,verbatim,amssymb}
\usepackage{times}       

 \newcommand{\resp}{{\it resp.} }
\newcommand{\cf}{{\it cf.} }
\newcommand{\ie}{{\it i.e.} }
\newcommand{\eg}{{\it e.g.} }

\newcommand{\Q}{\mathbb{Q}}
\newcommand{\R}{\mathbb{R}}
\newcommand{\C}{\mathbb{C}}
  \newcommand{\Z}{\mathbb{Z}}

\newcommand{\inj}{\hookrightarrow}

\newcommand{\car}{\operatorname{car}}
\newcounter{spec}

\swapnumbers

\newtheorem{thm}{Theorem}[subsection]

\newtheorem{prop}[thm]{Proposition}

\newtheorem{conj}[thm]{Conjecture}

\theoremstyle{definition}

\newtheorem{ex}[thm]{Example}
\newtheorem{exs}[thm]{Examples}

\newtheorem{rem}[thm]{Remark}

\numberwithin{equation}{section}

 at10pt
 at12pt

\begin{document}

\title[An introduction to motivic zeta functions of motives]{An introduction to motivic zeta functions of motives.
 } 
\author{Yves
Andr\'e}
 
  \address{D\'epartement de Math\'ematiques et Applications, \'Ecole Normale Sup\'erieure  \\ 
45 rue d'Ulm,  75230
  Paris Cedex 05\\France.}
\email{andre@dma.ens.fr}
  \date{\today}
  \begin{abstract}   It oftens occurs that Taylor coefficients of (dimensionally regularized) Feynman amplitudes $I$ with rational parameters, expanded at an integral dimension $D= D_0$, are not only periods (Belkale, Brosnan, Bogner, Weinzierl) but actually multiple zeta values (Broadhurst, Kreimer).

 In order to determine, at least heuristically, whether this is the case in concrete instances, the philosophy of motives - more specifically, the theory of mixed Tate motives - suggests an arithmetic approach (Kontsevich): counting points of algebraic varieties related to $I$ modulo sufficiently many primes $p$ and checking that the number of points varies polynomially in $p$. 
 
On the other hand, Kapranov has introduced a new ``zeta function", the role of which is precisely to ``interpolate" between zeta functions of reductions modulo different primes $p$. 

In this survey, we outline this circle of ideas and some of their recent developments.  

  \end{abstract}
\maketitle

 \renewcommand{\abstractname}{Summary}

 \begin{sloppypar}

  \bigskip    \bigskip

  This article is divided in two parts. 
  
  In the second and main part, we survey motivic zeta functions of motives, which ``interpolate" between Hasse-Weil zeta functions of reductions modulo different primes $p$ of varieties defined by polynomial equations with rational coefficients. 
  
  In the first and introductory part, we give some hints about the relevance of the concepts of motives and motivic zeta functions in questions related to computations of Feynman integrals. 
  

  \smallskip

 \bigskip
\section{Periods and motives.}

\subsection{Introduction.}  Relations between Feynman integrals and (Grothendieck) motives are manifold and mysterious. The most direct conceptual bridge relies on the notion of {\it period}, in the sense of arithmetic geometry; that is, integrals where both the integrand and the domain are defined in terms of polynoms with rational (or algebraic) coefficients.  The ubiquitous multiple zeta values encountered in the computation of Feynman amplitudes are emblematic examples of periods. 
  
 Periods are just complex numbers, but they carry a rich hidden structure reflecting their geometric origin. They occur as the entries of a canonical isomorphism between the complexification of two rational cohomologies attached to algebraic varieties defined over $\Q$: {\it algebraic De Rham cohomology}, defined in terms of algebraic differential forms\footnote{first defined by Grothendieck}, and {\it Betti\footnote{or singular} cohomology}, defined in terms of topological cochains. 
 
  Periods are thus best understood in the framework of {\it motives}, which are supposed to play the role of pieces of universal cohomology of algebraic varieties. For instance, the motive of the projective space ${\bf P}^n$ splits into $n+1$ pieces (so-called Tate motives) whose periods are $1, 2\pi i, \ldots, (2\pi i)^n$.
  
 \smallskip 
  To each motive over $\Q$ is associated a square matrix of periods (well defined up to left or right multiplication by matrices with rational coefficients), and a deep conjecture of Grothendieck predicts that the period matrix actually determines the motive. 
  
  For instance, multiple zeta values are periods of so-called mixed Tate motives over $\Z$,  which are iterated extensions of Tate motives. Grothendieck's conjecture implies that there are the only motives (over $\Q$) with multiple zeta values as periods.   
  
 \medskip In the philosophy of motives, cohomologies are thought of interchangeable realizations (functors with vector values), and one should take advantage of switching from one cohomology to another. Aside de Rham or Betti cohomology, one may also consider etale cohomology, together with the action of the absolute Galois group $Gal(\bar \Q/\Q)$; this amounts, more or less, to considering the number $N_p$ of points of the reduction modulo $p$ for almost all prime numbers $p$. 
 
A deep conjecture of Tate, in the same vein as Grothendieck's conjecture, predicts that the numbers $N_p$ determine the motive, up to semi-simplification.  For mixed Tate motives\footnote{not necessarily over $\Z$}, the $N_p$ are polynomials in $p$, and Tate's conjecture implies the converse. 

\medskip  To decide whether periods of a specific algebraic variety over $\Q$, say a hypersurface, are multiple zeta values may be a difficult problem about concrete integrals. The philosophy of motives suggests, as a test, to look at number of points $N_p$ of this hypersurface modulo $p$. Recently,  various efficient algorithms have been devised for computing $N_p$, \cf \eg \cite{Ke}\cite{St}. If $N_p$ turns out not to be polynomial in $p$, there is no chance that the periods are multiple zeta values (this would contradicts Grothendieck's or Tate's conjecture).

\subsection{Periods.}   A {\it period} is a complex number whose real and imaginary parts are absolutely convergent multiple integrals 
   $$\alpha = \frac{1}{{ \pi  }^m}\, \int_\Sigma \,\Omega$$ where $\Sigma $ is a domain in $\R^n$ defined by polynomial inequations with rational coefficients, $\Omega$ is a rational differential $n$-form with rational coefficients, and $m$ is a natural integer\footnote{the very name ``period" comes from the case of elliptic periods (in the case of an elliptic curve defined over $\bar\Q$, the periods of elliptic functions in the classical sense are indeed periods in the above sense) 
}. The set of periods is a countable subring of $\C$ which contains $\bar\Q$.
   
   \smallskip This is the definition proposed in \cite{KZ}\footnote{except for the factor $\frac{1}{{ \pi  }^m}$. We prefer to call {\it effective period} an integral  $\alpha$ in which $m=0$, to parallel the distinction motive versus effective motive }. There are some variants, which turn out to be equivalent. For instance, one could replace everywhere ``rational" by ``algebraic". Also one could consider a (non necessarily closed) rational $k$-form in $n$ variables with rational coefficients, integrated over a domain in $\R^n$ defined by polynomial equations and inequations with rational coefficients (interpreting the integral of a function as the volume under the graph, one can also reduce to the case when $k=n$ and $\Omega$ is a volume form).
   
    More geometrically, the ring of periods is generated by $\frac{1}{{ \pi }}$ and the numbers of the form 
    $\int_\gamma \,\omega$ where $\omega \in \Omega^n(X)$ is a top degree differential form on a smooth algebraic variety $X$ defined over $\Q$, and $\gamma\in H_n(X(\C),Y(\C); \Q)$ for some divisor $Y\subset X$ with normal crossings \cf \cite[p. 3, 31]{KZ}).
    
In many examples, \eg multiple zeta values (see below), the integrals have singularities along the boundary, hence are not immediately periods in the above sense. However, it turns out that any convergent integral $\int_\gamma \,\omega$ where $\omega \in \Omega^n(X\setminus Y)$ is a top degree differential form on the complement of a closed (possibly reducible) subvariety $Y$ of a smooth algebraic variety $X$ defined over $\Q$, and $\gamma$ is a semialgebraic subset of $X(\R)$ defined over $\Q$ (with non-empty interior), is a period
      (\cf \cite[th. 2.6]{BB2}).

   \subsection{Periods and motives.} Periods arise as entries of a matrix of the comparison isomorphism, 
 given by integration of algebraic diffential forms over chains, between 
 algebraic De Rham and ordinary Betti relative cohomology
   \begin{equation}\label{isocomp} H_{DR}(X,Y)\otimes \C \stackrel{\varpi_{X,Y}}{\cong}  H_{B}(X,Y)\otimes \C.\end{equation}
 $X$ being a smooth algebraic variety over $\Q$, and $Y$ being a closed (possibly reducible) subvariety\footnote{by the same trick as above, or using the Lefschetz's hyperplane theorem, one can express a period of a closed form of any degree as a period of a top degree differential form}.  
 
  \smallskip This is where motives enter the stage.  They are intermediate between algebraic varieties and their linear invariants (cohomology).    One expects the existence of an abelian category ${\rm MM(\Q)} $ of {\it mixed motives} (over $\Q$, with rational coefficients), and of a functor
    $$h:  Var(\Q) \to {\rm MM(\Q)}$$ (from the category of algebraic varieties over $\Q$) which plays the role of universal cohomology (more generally, to any pair $(X,Y)$ consisting of a smooth algebraic variety and a closed subvariety, one can attach a motive $h(X,Y)$ which plays the role of the universal relative cohomology of the pair).  
        
    The morphisms in ${\rm MM(\Q)}$ should be related to algebraic correspondences.
     In addition, the cartesian product on $Var(\Q)$ corresponds via $h$ to a certain tensor product $\otimes$ on ${\rm MM(\Q)}$, which makes ${\rm MM(\Q)}$ into a {\it tannakian category}, {\it i.e.} it has the same formal properties as the category of representations of a group. The (positive or negative)
 $\otimes$-powers  of $h^2({\bf P}^1)$ (and their direct sums) are called the pure {\it Tate motives}.
  
   \smallskip   The cohomologies $H_{DR}$ and $H_B$ factor through $h$, giving rise to two $\otimes$-functors 
 $$ H_{DR}, \, H_B:\, {\rm MM(\Q)} \to Vec_\Q$$ with values in the category of finite-dimensional $\Q$-vector spaces. 
   Moreover, corresponding to \eqref{isocomp},  there is an isomorphism in $Vec_\C$
 \begin{equation}  \varpi_M:  H_{DR}(M)\otimes \C \cong H_B(M)\otimes \C  \end{equation} which is
 $\otimes$-functorial in the motive $M$. The entries of a matrix of $\varpi_{M}$ with respect to some basis of the $\Q$-vector space $H_{DR} (M)$ (\resp $H_{B} (M )$) are the {\it periods} of $M$.
 
 \smallskip One can also consider, for each prime number $\ell$, the $\ell$-adic etale realization
 $$ H_{\ell} :\, {\rm MM(\Q)} \to Rep_{\Q_\ell}(Gal(\bar\Q /\Q) $$  with values in the category of finite-dimensional $\Q_\ell$-vector spaces endowed with continuous action of the absolute Galois group of $\Q$. As for varieties over $\Q$, it makes sense to reduce any motive $M\in MM(\Q)$ modulo sufficiently large primes $p$, and to ``count the number of points of $M$ modulo $p$". This number can be evaluated using the trace of Frobenius element at $p$ acting on $H_{\ell}(M)$.
 
       \subsection{Motivic Galois groups, period torsors, and Grothendieck's period conjecture.}  
  Let $\langle M \rangle$ be the tannakian subcategory of ${\rm MM(\Q)}$ generated by a motive $M$: its objets are given by algebraic constructions on $M$ (sums, subquotients, duals, tensor products).
  
  One defines the {\it motivic Galois group} of $M$ to be the group scheme
  \begin{equation} G_{mot}(M) := Aut^\otimes \,{H_B}_{\mid \langle M \rangle}  \end{equation}
  of automorphisms of the restriction of the $\otimes$-functor ${H_B}$ to $  \langle M \rangle$. 
  
This is a linear algebraic group over $\Q$: in heuristic terms, $ G_{mot}(M)$  is just {\it the Zariski-closed subgroup of $GL ( H_B(M))$  consisting of matrices which preserve motivic relations in the algebraic constructions on $H_B(M)$}.
  
  \smallskip Similarly, one can consider both $H_{DR}$ and $H_B$, and define the {\it period torsor} of $M$ to be 
      \begin{equation} P_{mot}(M) := Isom^\otimes \,({H_{DR}}_{\mid \langle M \rangle}, {H_B}_{\mid \langle M \rangle}) \in Var(\Q)  \end{equation}
  of isomorphisms of the restrictions of the $\otimes$-functors ${H_{DR}}$ and ${H_B}$ to $  \langle M \rangle$. 
This is a torsor under  $ G_{mot}(M) $, and it has a canonical complex point:
   \begin{equation}\varpi_M\in  P_{mot}(M)(\C).   \end{equation}

 \medskip{Grothendieck's period conjecture} asserts that  {\it the smallest algebraic subvariety of $P_{mot}(M)$ defined over $\Q$ and containing $\varpi$ is $P_{mot}(M)$ itself. }
  
  \smallskip 
 In more heuristic terms, this means that any polynomial relations with rational coefficients between periods should be of motivic origin (the relations of motivic origin being precisely those which define     $P_{mot}(M)$). This implies that a motive $M \in MM(\Q)$ can be recovered from its periods.

    \smallskip The conjecture is also equivalent to: 
 {\it  $P_{mot}(M)$ is connected (over $\Q$) and   
  \begin{equation}\label{trdegmot}{\rm tr. \,deg}_\Q \,  \Q[{\rm periods}(M)] = \dim \,G_{mot}(M) . \end{equation} }
  
  For further discussion, see \cite[ch. 23]{A1}.
 
\begin{ex} the motive of ${\mathbb P}^n$ decomposes as 
  \begin{equation} h({\mathbb P}^n)= \Q(0) \oplus \ldots \oplus \Q(-n), \end{equation} with periods $ 1, 2\pi i, \ldots, (2\pi i)^n$.
  Its motivic Galois group is the mutiplicative group ${\mathbb G}_m$. In this case, Grothendieck's conjecture amounts to the transcendence of $\pi$. \end{ex}
 
 \begin{rem}  By definition, periods are convergent integrals of a certain type. They can be transformed by algebraic changes of variable, or using additivity of the integral, or using Stokes formula. 
  
  M. Kontsevich conjectured that {\it any polynomial relation with rational coefficients between periods can be obtained by way of these elementary operations from calculus} (\cf \cite{KZ}). Using ideas of M. Nori \footnote{and granting the expected equivalence of various motivic settings}, it can be shown that this conjecture is actually equivalent to Grothendieck's conjecture (\cf \cite[ch. 23]{A1}).  \end{rem}
 
 Grothendieck's conjecture can be developped further into a Galois theory for periods, \cf \cite{A2}\cite{A3}.

\medskip
\subsection{Periods and Feynman amplitudes.} Let $\Gamma$ be a finite graph (without self-loop), with set of vertices $V$ and set of edges $E$. Let $\Psi_\Gamma$ be its classical Kirchhoff polynomial, \ie the homogeneous polynomial of degree $b_1(\Gamma)$ defined by  
\begin{equation}\Psi_\Gamma = \sum_T \prod_{e\notin T}\, x_e  ,\end{equation}
 where $T$ runs through the spanning trees of a given graph $\Gamma$ ($(x_e)$ is a set of indeterminates indexed by the edges of $\Gamma$).

Let $D_0$ be an even integer (for instance $4$). 
The graph $\Gamma$ can be considered as a (scalar) Feynman graph without external momenta. 
   According to the Feynman rules, when all masses are equal to $1$, the corresponding $D_0$-dimensional Feynman amplitude is written as
 \begin{equation} I_\Gamma(D_0) = \int_{\R^{D_0\vert E\vert}} \prod_{e\in E} (1+\vert p_e\vert^2)^{-1} \prod_{v\in V}\delta(\sum_{e\to v} p_e-\sum_{v\to e} p_e) \prod_{e\in E} d^{D_0}p_e.\end{equation}
 Its dimensional-regularization,  for $D$ close to $D_0$, can be evaluated, using the technique of Feynman parameters, to be 
 \begin{equation}I_\Gamma(D) =   \frac{\pi^{b_1(\Gamma)\cdot D/2}\cdot \Gamma(\vert E\vert -b_1(\Gamma)D/2)}{\Gamma(\vert E\vert)}\cdot J_\Gamma(D)\end{equation}
  where \begin{equation}J_\Gamma(D)= \int_{\Delta_{\vert E \vert}} \Psi_\Gamma^{-D/2}\prod_{e\in E} dx_e, \end{equation}  a convergent integral over the standard simplex $\Delta_{\vert E\vert}$ in $\vert E\vert$ variables. 

\smallskip In \cite{BB2}, P. Belkale and P. Brosnan showed that the Taylor coefficients of $J_\Gamma(D)$ at $D_0$ are periods (this is clear for $J_\Gamma(D_0)$ itself, except that one has to cope with the singularities of $\Psi_\Gamma^{-D/2}\prod_{e\in E} dx_e$ on $\Delta_{\vert E\vert}$, which requires a sequence of blow-ups. This issue is analysed in detail in \cite{BEK}\cite{BW1}\cite{M}. For other Taylor coefficients, one has to add one more variable, \cf \cite[p. 2660]{BB2}).
 
\smallskip In \cite{BW}, this result was extended to the case when $\Gamma$ is a semi-graph (\ie in the presence of external momenta) and when the masses are non necessarily equal to $1$ but are commensurable to each other. 
  
\begin{rem} Taking into account these results, polynomial relations between Feynman amplitudes attached to different graphs $\Gamma$ (like the relations which lead to Kreimer's Hopf algebra) can be intepreted as period relations. According to Grothendieck's conjecture, they should be of motivic origin, \ie come from relations between the motives attached to the hypersurfaces $X_\Gamma: \Psi_\Gamma=0$ (and related varieties). Some evidence for this is given in \cite{BEK}.\end{rem}

\medskip
\subsection{Multiple zeta values, Feynman amplitudes and Hasse-Weil zeta functions.} 

 Multiple zeta values 
\begin{equation}\label{polyzetadef} \zeta(s_1,\dots, s_k)=\sum_{n_1>\dots
>n_k\ge 1} \,\frac{1}{n_1^{s_1}\dots n_k^{s_k}}\end{equation}
(where $s_i$ are integers $\ge 1$, with $s_1\ge 2$) can be written in integral form:
 setting
\[ \omega_0 = \frac{dt}{t}\;\;  \omega_1 = \frac{dt}{1-t},\;\;\omega_r=
\omega_0^{\wedge (r-1)}\wedge\omega_1 \;\text{pour }\; r\ge 2,\]   one has
\begin{equation}\label{polylogint}\zeta(s_1,\dots, s_k)= \int_{1>t_1>\dots >t_{s}>0} \omega_{s_1} \dots
\omega_{s_k},\end{equation} 
which is thus a period. This is actually the period of a mixed Tate motive over $\Z$, \ie an iterated extension in $MM(\Q)$ of pure Tate motives, which is unramified with respect to the Galois action on etale cohomology (\cf \cite[ch. 25]{A1} for more detail). 

\smallskip These numbers have long been known to occur in a pervasive manner in the computation of Feynman amplitudes (\cf \eg \cite{BK}, \cite{W}). 

\smallskip Kontsevich once speculated that the period of the hypersurface $X_\Gamma: \Psi_\Gamma= 0$ were 
  (linear combinations of) multiple zeta values.  According to Grothendieck's period conjecture, this would imply that the motive of $X_\Gamma$ is a mixed Tate motive over $\Z$.  
 If this is the case, the number $\sharp X_\Gamma(\mathbb F_p)$ of points of the reduction of $X_\Gamma$ modulo $p$ should be polynomial in $p$; equivalently, the poles of the Hasse-Weil zeta function of $X_\Gamma\otimes \mathbb F_p$ should be integral powers of $p$. 

This has been checked for graphs with less than $12$ edges by J. Stembridge \cite{St}, but disproved in general by Belkale and Brosnan \cite{BB1}. 

\medskip
 However, this leaves open the interesting general question, for any $X\in Var(\Q)$, of {\it controlling $\sharp X(\mathbb F_p)$ uniformly in $p$ - or equivalently, of the variation of $Z(X\otimes \mathbb F_p, t)$ with $p$}. 
 
As we shall see, there are well-suited mathematical tools to tackle this question:  the {\it Kapranov zeta function}, and its variant the {\it motivic zeta function}. 

 \begin{rem} The relationship between Feynman diagrams  and motives has been investigated much further. For instance, in \cite{AM}, P. Aluffi and M. Marcolli propose an algebro-geometric version of the Feynman rules, which takes place in a certain $K_0$-ring built from immersed conical varieties.  
\end{rem}

 \medskip
 
 \bigskip\section{Motivic zeta functions.}
 
\medskip\subsection{The ring of varieties.} The idea to build a ring out of varieties, viewing pasting as the addition, is very old. In the case of algebraic varieties over a field $k$, this leads to the ring $K_0(Var(k))$\footnote{which occurs in some early letters from Grothendieck to Serre about motives}: the generators are denoted by $[X]$, one for each isomorphism class of $k$-variety; the relations are generated by 
\begin{equation}[X-Y]= [X]-[Y]\end{equation}
when $Y$ is a closed subvariety of $X$. With the product given by 
\begin{equation}[X\times Y]= [X]\cdot [Y],\end{equation}
$K_0(Var(k))$ becomes a ring. 

It is standard to denote by $\mathbb L$ the class $[{\bf A}^1]$ of the affine line.

\begin{exs}: 1) One has 
\begin{equation} [GL_n] = ({\mathbb L}^n-1)\ldots ({\mathbb L}^n-{\mathbb L}^{n-1})=  (\mathbb L - 1)\cdot [SL_n].\end{equation}

\noindent 2) In the case of a Zariski locally trivial fibration $X\to S$ with fiber $Y$, one has $[X]=[S]\cdot[Y].$
This applies to $GL_n$-fibrations (which are locally trivial); in order to recover $[S]$ from $[X]$, taking into account the formula for $[GL_n]$, it will often be convenient to localize $K_0(Var(k))$ by $\mathbb L$ and $\mathbb L^n- 1,\, n>0$. 

\noindent 3) (McWilliams, Belkale-Brosnan \cite{BB1}): the class of the space of forms of rank $r$ in $n$ variables is
\begin{equation}[Sym^n_r] = \prod_1^s \frac{\mathbb L^{2i}}{\mathbb L^{2i}-1 }\prod_0^{2s-1}(\mathbb L^{n-i} -1)\end{equation}
if $0\leq r=2s\leq n$,

\begin{equation}[Sym^n_r] = \prod_1^s \frac{\mathbb L^{2i}}{\mathbb L^{2i}-1 }\prod_0^{2s }(\mathbb L^{n-i} -1)\end{equation}
if $0\leq r=2s+1\leq n$.
\end{exs}

The structure of $K_0(Var(k))$ is rather mysterious. It is slightly better understood when $k$ is of characteristic $0$, using strong versions of the resolution of singularities. 

\begin{prop}\cite{Bi} If $\car k=0$, $K_0(Var(k))$ admits the following presentation: generators are classes of smooth projective varieties $X$, with the blow-up relations:
\begin{equation}[Bl_Y X] -[E]= [X]-[Y] \end{equation}
where $E$ denotes the exceptional divisor in the blow-up $Bl_Y X$ of $X$ along the smooth subvariety $Y$.\end{prop}
  
\medskip\subsection{Relation to motives.}  In the category $MM(k)$ of mixed motives over $k$ with rational coefficients\footnote{there are actually several candidates for this category, some conditional, some not}, relations (2.1), (2.2), (2.3), (2.6) have more sophisticated counterparts, which actually reduce to analogous relations if one passes to $K_0(MM(k)$, the Grothendieck group constructed in terms of  extensions of mixed motives.  

In fact, one expects that the functor\footnote{more accurately, its variant with compact supports.  } $h: Var(k) \to MM(k)$ gives rise to a ring homomorphism 
\begin{equation} K_0(Var(k)) \to K_0(MM(k)).\end{equation}
This can be made rigorous, if $\car k=0$, using the previous proposition and a category $M_\sim(k)$ of {\it pure motives} (\ie of motives of smooth projective $k$-varieties, with morphisms given by algebraic correspondences modulo some fixed equivalence relation $\sim$). One gets a canonical ring homomorphism 
\begin{equation} \mu_c:\, K_0(Var(k)) \to K_0(M_\sim(k)) \end{equation}
(where $K_0(M_\sim(k))$ denotes the Grothendieck group constructed in terms of direct sums of pure motives; which is actually a ring with respect to the multiplication induced by tensoring motives). It sends $\mathbb L$ to $[\Q(-1)]$.

Recent work by H. Gillet and C. Soul\'e \cite{GS} allows to drop the assumption on $\car k$.

\begin{rem}
Conjecturally, $K_0(M_\sim(k))= K_0(MM(k))$ and does not depend on the chosen equivalence $\sim$  used in the definition of the $M_\sim(k)$. In fact, this independence would follow from a conjecture due to S. Kimura and P. O'Sullivan, which predicts that any pure motive $M\in M_\sim(k)$ decomposes (non-canonically) as $M_+\oplus M_-$, where $S^n M_-= \bigwedge^n M_+ = 0$ for $n>>0$ (\cf \eg \cite[ch. 12]{A1}\footnote{for the coarsest equivalence $\sim$ (the so-called numerical equivalence), this conjecture amounts to the following:
the even K\"unneth projector  $ H(X)\to H^{even}(X)\to H(X)$ is algebraic}; here $S^n$ and $\bigwedge^n$ denote the $n$-th symmetric and antisymmetric powers, respectively). This is for instance the case for motives of products of curves.

In the sequel, we shall deal with sub-$\otimes$-categories of $M_\sim(k)$ which satisfy this conjecture, and we will drop $\sim$ from the notation $K_0(M_\sim(k))$. 
\end{rem}

\medskip\subsection{Kapranov zeta functions.}  
When $k$ is a finite field, counting $k$-points of varieties factors through a ring homomorphism
\begin{equation} \nu :\, K_0(Var(k)) \to \Z,\;\; [X]\mapsto \sharp X(k),\end{equation}
which factors through $K_0(M(k))$. One of the expressions of the {\it Hasse-Weil zeta function}, which encodes the number of points of $X$ in all finite extensions of $k$, is
\begin{equation}Z(X,t)=\sum_0^\infty\,  \sharp ((S^nX)(k)) \, t^n \in \mathbb Z[[t]]  ,\end{equation}
 and it belongs to $\Q(t)$ (Dwork).

\medskip
M. Kapranov had the idea \cite{Kap} to replace, in this expression, $\sharp ((S^nX)(k))$ by the class of $S^nX$ itself in $K_0(Var(k))$\footnote{for $X$ quasiprojective, say, in order to avoid difficulties with symmetric powers}). More precisely, he attached to any ring homomorphism
$$\mu: K_0(Var(k) \to R  $$ the series 
\begin{equation} Z_\mu(X,t):= \sum_0^\infty\,\mu[S^nX] \, t^n \in R[[t]],\end{equation}
 which satisfies the equation
  $$ Z_\mu(X \coprod X' ,t) = Z_\mu(X   ,t).Z_\mu( X' ,t) .$$

When $\mu= \nu$ (counting $k$-points), one recovers the Hasse-Weil zeta function. 

When  $k=\mathbb C$ and $ \mu = \chi_c $ (Euler characteristic), $ Z_\mu(X,t) = (1-t)^{-\chi_c(X)}$ (MacDonald). 

The universal case (the {\it Kapranov zeta function}) corresponds to $\mu= id$. When $k= \Q$, one can reduce $X$ modulo $p>>0$ and count $\mathbb F_p$-points of the reduction. The Kapranov zeta function then specializes to the Hasse-Weil zeta function of the reduction, and thus may be seen as some kind of interpolation of these Hasse-Weil zeta functions when $p$ varies.

\medskip\subsection{Around the Kapranov zeta function of a curve.} 
Let us assume that $X$ is a smooth projective curve of genus $g$, defined over the field $k$. The Kapranov zeta function 
\begin{equation} Z_\mu(X,t):= \sum_0^\infty\,\mu[S^nX] \, t^n\end{equation}
has the same features as the usual Hasse-Weil zeta function:
 
\begin{prop}\cite{Kap} \begin{equation} Z_\mu(X,t)= \frac{P_\mu( X, t)}{(1-t)(1-\mathbb L t)}\end{equation}
 where $P_\mu$ is a polynomial of degree $2g $, and 
one has the functional equation  \begin{equation}Z_\mu(X,t)=\mathbb L^{g-1} t^{2g-2} Z_\mu(X,\mathbb L^{-1} t^{-1}). \end{equation}

     \end{prop}
 
 Sketch of proof of (2.13): the mapping
$$X^n \to J(X),\; (x_1,\ldots, x_n)\mapsto [x_1]+\ldots [x_n]-n[x_0]$$
factors through $ S^nX \to J(X), $ which is a projective bundle if $n\geq 2g-1$.

Moreover, one has an injection
$$S^n(X)\inj S^{n+1}(X),\;\; x_1+\ldots x_n\mapsto x_0+x_1+\ldots x_n$$
and the complement of its image is a vector budle of rank $n+1-g$ over $J(X)$.
 This implies
  $ [S^{n+1}(X)]-[S^n(X)]= [J(X)]\mathbb L^{n+1-g}$ 
 hence, by telescoping, that
  $Z_\mu(X,T)(1-T)(1-\mathbb L T)$  is a polynomial of degree $\leq 2g$. $\square$

When $k$ is a finite field, the Hasse-Weil zeta function of $X$ can also be written in the form
  \begin{equation} \displaystyle Z(X,t) = \sum_{   D\geq 0}\, t^{\deg D}  \;\; \text{\rm(sum over effective divisors)}\end{equation}
    \begin{equation} = \sum_{\mathcal L} \, h^0(\mathcal L)_q .\, t^{\deg \mathcal L}  \;\; \text{\rm(sum over line bundles}),\end{equation}
 where one uses the standard notation   $n_q= 1+q+\ldots + q^{n-1}$.
 
 R. Pellikan \cite{P} had the idea to substitute, in this expression, $q$ by an indeterminate $u$. He proved that  
  \begin{equation} Z(X,t,u) :=  \sum_{\mathcal L} \, h^0(\mathcal L)_u .\, t^{\deg \mathcal L}\end{equation}
  is a rational function of the form $\;\frac{P(X, t, u)}{(1-t)(1-u t)}$.
 
  Finally, F. Baldassarri, C. Deninger and N. Naumann \cite{BDN} unified the two generalizations (2.12) (Kapranov) and (2.17) (Pellikaan) of the Hasse-Weil zeta function (2.15) by setting:
 \begin{equation} Z_\mu(X,t,u) := \sum_{n,d} \, [Pic^d_n] n_u .\, t^{d} \in R[[t,u]]\end{equation} (where $Pic^d_n$ classifies line bundles of degree $d$ with $h^0(\mathcal L)\geq n$, and $n_u= 1+u+\ldots + u^{n-1}$),
  and they proved that this is again a rational function of the form $\;\frac{P_\mu(X, t, u)}{(1-t)(1-u t)}$.
 
 One thus has a commutative diagram of specializations
   $$\begin{matrix} &&Z_{\mu}(X,t)&& \\ &\stackrel{u\mapsto q}{\nearrow}  &  &  \stackrel{\mu\mapsto\nu}{\searrow} &\\ Z_{\mu}(X,t,u) && &&  Z(X,t). \\  &\stackrel{\mu\mapsto\nu}{\searrow}  &  &  \stackrel{u\mapsto q}{\nearrow} &\\ &&Z (X,t,u)&& \end{matrix}$$
  
\smallskip On the other hand, M. Larsen and V. Lunts investigated the Kapranov zeta function of products of curves.

\begin{prop} \cite{LL1}\cite{LL2} If $X$ is a product of two curves of genus $>1$, 
 $Z_{\mu}(X,t) $ is \emph{not} rational for $\mu=id$. \end{prop}
 
In the sequel, following \cite[13.3]{A1}, we remedy this by working with a $\mu$ which is ``sufficiently universal", but for which one can nevertheless hope that $Z_\mu(X,t)$ is always rational. Namely, we work with $\mu_c$: in other words, we replace the ring of varieties by the $K$-ring of pure motives.

\medskip\subsection{Motivic zeta functions of motives.}  
    Thus, let us define, for any pure motive $M$ over $k$ (with rational coefficients), its {\it motivic zeta function} to be the series
   \begin{equation} Z_{mot}(M,t):= \sum_0^\infty\, [S^n M] . \, t^n  \in K_0({M(k)} )[[t]].\end{equation}
  One has 
 $$ Z_{mot}(M \oplus M' ,t) = Z_{mot}(M   ,t).Z_{mot}( M' ,t) .$$

\begin{prop}\cite[13.3]{A1} If $M$ is finite-dimensional in the sense of Kimura-O'Sullivan (\ie $M= M_+\oplus M_-$,\; $S^{n}M_-= \bigwedge^n M_+= 0$ for $n>>0$), then $Z_{mot}(M,t)$ is rational.\end{prop}
  
(This applies for instance to motives of products of curves - and conjecturally to any motive). 
 
\smallskip Moreover, B. Kahn \cite{Kah} (\cf also \cite{He}) has established a functional equation of the form 
\begin{equation} Z_{mot}(M^\vee   ,t^{-1})  =   (-1)^{\chi_+(M)}. \det M. t^{\chi(M)} . Z_{mot}(M  , t)\end{equation}
 (where $\det M = \bigwedge^{\chi_+}M_+ \otimes (S^{-\chi_-}M_-)^{-1}$).

\bigskip\subsection{Motivic Artin $L$-functions.} One can play this game further. Hasse-Weil zeta functions of curves can be decomposed into (Artin) $L$-functions. A. Dhillon and J. Min\'ac upgraded this formalism at the level of motivic zeta functions \cite{DM}. 

\smallskip Starting in slightly greater generality, let $V$ be a $\Q$-vector space of  finite dimension. To any pure motive $M$, one attaches another one $V\otimes M$, defined by 
  $$Hom(V\otimes M, M') = Hom_F(V, Hom(M,M')).$$
Let $G$ be a finite group, and let $\rho:\, G \to GL(V)$ be a homomorphism. 

The {\it motivic $L$-function} attached to $M$ and $\rho$ is 
   \begin{equation} L_{mot}(M,\rho,  t) := Z_{mot}((V\otimes M)^G, t) .\end{equation}

This definition extends to characters $\chi$ of $G$ (that is, $\mathbb  Z$-linear combinations of $\rho$'s), and gives rise to the a formalism analogous to the usual formalism of Artin $L$-functions. Namely, one has the following identities in $K_0(M(k))(t)$:
      \begin{equation}L_{mot}(M, \chi+\chi', t) = L_{mot}(M,\chi, t).L_{mot}(M,\chi', t)   \end{equation}
   \begin{equation}L_{mot}(M,  \chi', t) = L_{mot}(M,Ind_{G'}^G\chi', t) \;   \end{equation}  (for $G'$ a subgroup of $G$),
   \begin{equation}L_{mot}(M,  \chi'', t ) = L_{mot}(M,\chi,t) \;    \end{equation}  ($G'\triangleleft G$, $\chi$ coming from a character $\chi''$ of $G/G'$),
   \begin{equation}Z_{mot}(M,   t) = \prod_{ \chi \,{\rm{irr.}}} \, L_{mot}(M,\chi, t)^{\chi(1)}.   \end{equation}

\begin{ex} Let $X$ be again a smooth projective curve, and let $G$ act on $X$. Then $G$ acts on the motive $h(X)$ of $X$ via $(g^\ast)^{-1}$. By definition $L_{mot}(X,\chi, t)$ is the motivic $L$-function of $h(X)$.
 
 If $k$ is finite, $\nu(L_{mot}(X,\chi, t))$ is nothing but the  Artin non-abelian $L$-function $L(X, \chi, t)$, defined by the formula (where  $F$ denotes the Frobenius)
$$\log L(X, \chi, t)= \sum \nu_n(X) \frac{t^n}{n}\;\nu,\;\; \nu_n(X)= \frac{1}{\sharp G}\sum\chi(g^{-1})\sharp Fix(gF^n).$$ \end{ex}
 
This leads to the definition of {\it motivic Artin symbols} and to a motivic avatar of  Cebotarev's density theorem \cite{DM}.

\medskip \subsection{The class of the motive of a semisimple group $G$}
 Let $G$ be a connected split semisimple algebraic group over $k  $. Let
$T\subset G$ be a maximal torus (of dimension $r$), with character group $X(T)$.
  The Weyl group $W$ acts on the symmetric algebra $S(X(T)_\Q))$, and its invariants are generated in degree $d_1, \ldots, d_r $ (for $G=SL_n$, one has $ d_1= 2,\ldots, d_r=n$).

   One has the classical formulas

\begin{equation}(t-1)^r\sum_{w\in W} t^{\ell(w)}= \prod_1^r\, (t^{d_i}-1),\;\;
 \sum d_i = \frac{1}{2} (\dim G+r).\end{equation}

K. Behrend and A. Dhillon gave the following generalization of formula (2.3) for $[SL_n]$.

\begin{prop} In $K_0(Var(k))[\mathbb L^{-1}]$ or $K_0({M(k))}$, one has 
\begin{equation} [G]= \mathbb L^{\dim G}\prod_1^r (1-\mathbb L^{-d_i}).\end{equation}
\end{prop} 
(In $K_0({M(k))}$, it is preferable to write $[\Q(m)]$ instead of $\mathbb L^{-m}$).

Sketch of proof: let $B$ be a Borel subgroup, $U$ its unipotent radical. Then
$  [G]= [G/B]\cdot[T]\cdot[U] $ (in $K_0(Var(k))$ or $K_0({M(k))}$). On the other hand, one computes easily
$  [U]= \mathbb L^{ \frac{1}{2} (\dim G-r)},\; [T]= (\mathbb L-1)^r $, and using the Bruhat decomposition,
$\,   [G/B] = \sum_{w\in W} \mathbb L^{\ell (w)}.$  
Combining these formulas with (2.26), one gets (2.27). $\square$

\medskip With proper interpretation, (2.27) can be reformulated as a formula for the class of the classifying stack of $G$-torsors over $k$, in a suitable localization of $K_0(M(k))$:
\begin{equation}[BG]= [G]^{-1}= [\Q(\dim G)]\prod_i\,(1-[\Q(d_i)])^{-1}.\end{equation}

 \medskip\subsection{$G$-torsors over a curve $X$, and special values of $Z_{mot}(X)$}  
Let us now look at $G$-torsors not over the point, but over a smooth projective curve $X$ of genus $g$. 

More precisely, let $G$ be a simply connected split semisimple algebraic group over $k$ (\eg  $SL_n$), and let  $Bun_{G,X}$ be the moduli stack of $G$-torsors on $X$  (which is smooth of dimension $(g-1).\dim G$).

This stack admits a infinite stratification  by pieces of the form $[X_i/GL_{n_i}]$, whose dimensions tend to  $ -\infty$. According to Behrend and Dhillon \cite{BD}, this allows to define unambiguously the class

\begin{equation} [Bun_{G,X}]:= \sum [X_i][GL_{n_i}]^{-1}  \end{equation} in a suitable completion of $K_0[{M(k)}]$ with respect to $[\Q(1)]$, taking into account the fact that  $[GL_{n }]^{-1} = [BGL_{n}] =  \Q(n^2)\cdot( 1 +\cdots ) \in \mathbb Z [[  \Q(1) ]]$.

\begin{conj} (Behrend-Dhillon)   \begin{equation} [Bun_{G,X}]= [\Q((1-g).\dim G)]\prod_i Z_{mot}(X,[\Q(d_i)]).\end{equation} \end{conj}
\noindent This has to be compared with (2.28), where the $d_i$ have the same meaning;. Note that the special values $Z_{mot}(X, [\Q(d_i)])$ are well-defined since $Z_{mot}(X,t)$ is rational with poles at $1$ and $[\Q(1)]$ only. 

\begin{prop}\cite{BD} The conjecture holds for $X= \mathbb P^1$ and any $G$, and for $G=SL_n$ and any $X$.\end{prop}

Let us consider the case of $SL_n$ to fix ideas (\cf \cite{D}), and comment briefly on some specializations of the motivic formula (2.30).

\smallskip 1) For $k=\mathbb C,\, \mu = \chi_c$, the formula specializes to a formula for the Euler characteristic of $Bun_{G,X}$, which can be established via gauge theory \`a la Yang-Mills (Atiyah-Bott \cite{AB}, see also Teleman \cite{Te}). 

\smallskip More precisely, $H^*(Bun_{G,X} ) \cong H^*(G)^{\otimes 2g}\otimes H^*(BG)\otimes H^*(\Omega G)$.

\medskip 2) For $k=\mathbb F_q,\;\mu=\nu$ (counting points), the formula specializes to a formula for the number of $k$-points of $Bun_{G,X}$ (Harder, \cf \cite{HN}). 

\smallskip More precisely, $Bun_{G,X}$ can be viewed as the transformation groupoid of $G(K) $ on $G({\mathbb A}_K)/\mathcal K$, for  $K=k(X),\;\mathcal K = \prod_x\, G(\hat {\mathcal O}_{X,x})$; so that
 $\;\sharp Bun_{G,X} (k) =  \frac {vol(G(K)\backslash G({\mathbb A}_K))}{vol(\mathcal K)}$. One has 
 $\; \displaystyle vol (\mathcal K)  = q^{(1-g).(n^2-1)}\prod_2^n \zeta_K(i)^{-1},\,$ 
and the Tamagawa number
$\;\displaystyle vol(G(K)\backslash G({\mathbb A}_K)) $ is $1$,  whence
\begin{equation} \sharp Bun_{G,X} (k) = q^{(g-1).(n^2-1)}\prod_2^n \zeta_K(i). \end{equation} 
 
\end{sloppypar} 
\end{document}